 \documentclass[12pt,{other stuff},reqno,oneside]{amsart} 
  
\usepackage{amsmath,amscd}   
\usepackage{amssymb}

\usepackage[all]{xy}

\usepackage{amsthm}  
\newtheorem{theorem}{Theorem}[section]
\newtheorem{lemma}[theorem]{Lemma}
\newtheorem{proposition}[theorem]{Proposition}
\newtheorem{corollary}[theorem]{Corollary}
\newtheorem{definition}{Definition}[section]

\newtheorem{remark}{Remark}[section]
\def\dom{\textbf{dom}}

\def\r{\textbf{r}}
\def\d{\textbf{d}}
\def\Ob{\text{Ob\,}}
\def\Hom{\text{Hom\,}}
\def\cod{\textbf{cod}}

\usepackage{mathrsfs}

\usepackage[usenames,dvipsnames]{color}
   \usepackage[pdftex]{graphicx}
   
\usepackage{pb-diagram,lamsarrow,pb-lams}

  \usepackage{wrapfig}
 \usepackage{caption}{\upshape}
 
\begin{document}

\author{N. Ghroda}
  \address{Department of Mathematics\\University
  of York\\Heslington\\York YO10 5DD\\UK}\email{ng521@york.ac.uk}
 		
	 \title[Groupoids of left quotients]{Groupoids of left quotients}

 \begin{abstract}
     A subcategory $\textbf{C}$  of a groupoid  $\mathbb{G}$  is a left order in  $\mathbb{G}$,  if every element of  $\mathbb{G}$  can be written as  $a^{-1}b$    where  $a ,b \in \textbf{C}$. A subsemigroupoid $\mathfrak{C}$ of a groupoid $\mathbb{G}$ is a left q-order in  $\mathbb{G}$,  if every element of  $\mathbb{G}$  can be written as  $a^{-1}b$    where  $a ,b \in \mathfrak{C}$. We give a characterization of left orders (q-orders) in groupoids. In addition, we describe the relationship between left I-orders in primitive inverse semigroups and left orders (q-orders) in groupoids.
   \end{abstract}
   
 \keywords{primitive inverse semigroup, groupoids of left quotients, I-quotients, I-order}
 \date{\today}
%\large
%In the name of God

\maketitle 

\section{Introduction}

In this article we investigate left orders (q-orders) in groupoids. This work is part of a continuing investigation of categories of quotients. The motivation for our investigation comes from  semigroups of quotients and categories of fractions. Our purpose is the investigation of a similar problem in groupoid theory.\par
\medskip
 Fountain and Petrich introduced the notion of a completely 0-simple semigroup of quotients in \cite{pjhon}. It is well-known that groupoids are generalisations of groups, also, inverse semigroups can be regarded as special kinds of ordered groupoids. The concept of semigroups of quotients extends that of a group of quotients, introduced by Ore-Dubreil.  We recall that a group $G$ is a \emph{group of left quotients} of its subsemigroup $S$ if every element of $G$ can be written as $a^{-1}b$ for some $a,b\in S$. \par
\medskip
Ghroda and Gould \cite{GG2} extended the classical notion of left orders in inverse semigroups. They have introduced the following definition: Let  $Q$ be an inverse semigroup. A  subsemigroup  $S$  of  $Q$  is a \emph{left I-order} in  $Q$ and $Q$ is a \emph{semigroup of left I-quotient} of $S$, if every element of  $Q$  can be written as  $a^{-1}b$    where  $a ,b \in S$  and  $a^{-1}$  is the inverse of  $a$  in the sense of inverse semigroup theory. The notions of \emph{right I-order} and \emph{semigroup of right I-quotients} are defined dually. If  $S$ is both a left and a right I-order in an inverse semigroup  $Q$, we say that  $S$  is an \emph{I-order} in  $Q$ and $Q$ is a semigroup of \emph{I-quotients} of $S$. If we insist on  $a$ and  $b$  being  $\mathcal{R}$-related in  $Q$, then we say that  $S$ is a \emph{straight left I-order} in  $Q$. \par
 \medskip
 The theory of categories of fractions was developed by Gabriel and Zisman \cite{Gaberiel}.  The key idea is that starting with a category $\textbf{C}$ we can associate a groupoid to $\textbf{C}$ by adding all the inverses of all the elements of $\textbf{C}$ to $\textbf{C}$. We then produce a groupoid  $\mathbb{G}(\textbf{C})=\textbf{C}^{-1}\textbf{C}$ and a functor $\iota: \textbf{C} \longrightarrow \mathbb{G}$ such that $\mathbb{G}(\textbf{C})$ is generated by $\iota(\textbf{C})$, we call $\mathbb{G}(\textbf{C})$ a \emph{category of fractions}. Tobias in \cite{Tob} showed that for any category with conditions which are analogues of the Ore condition in the theory of non-commutative rings (see, \cite{Jat}), there is a groupoid of fractions.    \par
 \medskip
 Now, we are in a position to define a groupoid of left quotients. Let $\textbf{C}$ be a subcategory of a groupoid $\mathbb{G}$. We say that $\textbf{C}$ is a \emph{left order} in $\mathbb{G}$ or $\mathbb{G}$ is a groupoid of \emph{left quotients} of $\textbf{C}$ if  every element of $\mathbb{G}$ can be written as $a^{-1}b$ for some $a,b\in \textbf{C}$. \emph{Right orders} and \emph{groupoids of right quotients} are defined dually. If $\textbf{C}$ is both a left and a right order in $\mathbb{G}$, then $\textbf{C}$ is an \emph{order} in $\mathbb{G}$ and $\mathbb{G}$ is a \emph{groupoid of quotients} of $\textbf{C}$.\par
 \medskip 
In the case where the category without identities also known as a `\emph{quiver}' or a \emph{semigroupoid}, we propose the following definition of groupoid of q-quotients.\par
 \medskip
Let $\mathfrak{C}$ is a semigroupoid. We say that a groupoid $\mathbb{G}$ contains $\mathfrak{C}$ is a \emph{left q-qutients} of $\mathfrak{C}$ or $\mathfrak{C}$ is a \emph{left q-order} in $\mathbb{G}$ if  every element  of $\mathbb{G}$ can be written as $a^{-1}b$ for some $a,b\in \mathfrak{C}$. \emph{Right q-orders} and \emph{groupoids of right q-quotients} are defined dually. If $\mathfrak{C}$ is both a left and a right q-order in $\mathbb{G}$, then $\mathfrak{C}$ is a \emph{q-order} in $\mathbb{G}$ and $\mathbb{G}$ is a \emph{groupoid of q-quotients} of $\mathfrak{C}$.\par
\medskip
This work is divided up into eight sections. In Section~\ref{prel}  we summarize the background on groupoids and inverse semigroups that we shall need throughout the article. A Theorem 1.24 in \cite{clifford} due to Ore and Dubreil shows that a semigroup $S$ has a group of left quotients if and only if it is \emph{right reversible}, that is,  $Sa\cap Sb\neq \emptyset$ for all $a,b\in S$ and $S$ is cancellative. In Section~\ref{leftorder} we prove the category version of such a theorem. We stress that this work is not new - it has been studied by a number of authors,  by using the notion of category as a collection of objects and arrows. We regard a small category as a generealisation of a moinoid to prove such a theorem. Consequently, the relationship between the groupoids of left quotients and inverse semigroups of left I-quotients becomes clearer.\par
\medskip 
In Section~\ref{uniq}  we show that a groupoid of left quotients is  unique up to isomorphism. In Section~\ref{primitiveord} we establish the connection between groupoids of left quotients and primitive inverse semigroups of left I-quotients, then we specialise our result to left orders in connected groupoids in Section~\ref{connectedored}  a result which may be regarded as a generalisation of  Corollary 3.11 in \cite{GG}, which characterised left I-orders in Brandt semigroups. Recall that if we adjoin an element $0$ to a groupoid $\mathbb{G}$  and declare undefined products in $\mathbb{G}$  as being equal to zero, we obtain a primitive inverse semigroup $\mathbb{G}^0$.  We note that if a category $\textbf{C}$ is a left order in $\mathbb{G}$, then  $S=C^0$ is a left I-order  $Q=\mathbb{G}^0$. In fact, $S$ is a \emph{full subsemigroup} of $Q$ in the sense that $E(S)=E(Q)$. In some sense, this justifies introducing the notion of left q-order and left q-quotients.\par
 \medskip
 The structure of  inductive groupoids that correspond to 	 bisimple inverse $\omega$-semigroups was determined in \cite{Gel}. Then it is a natural question to ask for the relationship between  left orders in such  inductive groupoids and  left I-orders in bisimple inverse $\omega$-semigroups. In Section~\ref{w-ordered} we investigate such a relationship.\par
 \medskip
Theorem~\ref{cateqqu}, in Section~\ref{leftq-order} gives necessary and sufficient conditions for a semigroupoid to have a groupoid of left q-quotients.

\section{Preliminaries and notion}\label{prel}

 In this section we set up the definitions and results about groupoids and inverse semigroups. Standard references include \cite{clifford} for inverse semigroups, and \cite{Higgins} for groupoids. \par
 \medskip
There are two definitions of (small) category. The first one in \cite{Higgins} considers the category as a collection of objects (sets) and homomorphisms between them satisfying certain conditions. A \emph{category}, consisits of a set of objects $\{a,b,c, . . . \}$ and homomorphisms
between the objects such that: \\

$(i)$ Homomorphisms are composable: given homomorphisms $a: u \longrightarrow v$ and $b : v \longrightarrow w$, the homomorphism $ab: u \longrightarrow w$ exists, otherwise $ab$ is not defined;

$(ii)$ Composition is associative: given homomorphisms  $a: u \longrightarrow v$, $b : v \longrightarrow w$, and $c: w \longrightarrow z$, $(ab)c = a(bc)$;

$(iii)$ Existence of an identity homomorphism:  For each object $u$, there is an identity homomorphism $e_u: u \longrightarrow u$  such that for any homomorphism $a: u \longrightarrow v, \ e_ua = a= a e_v$.

\par
\medskip  
A category  $\mathbf{T}$ is called a  \emph{subcategory} of the category $\mathbf{C}$, if the objects of $\mathbf{T}$ are also objects of $\mathbf{C}$, and the homomorphisms of $\mathbf{T}$ are also homorphisms of $\mathbf{C}$ such that \\

    $(i)$ for every $u$ in $\Ob (\textbf{T})$, the identity homomorphism  $e_u$ is in $\Hom \textbf{T}$;
    
    $(ii)$ for every pair of homomorphisms $f$ and $g$ in $\Hom \textbf{T}$ the composite $fg$ is in $\Hom \textbf{T}$ whenever it is defined.\\
\\

The second definition regards the category as an algebraic structure in its own right. In this definition we can look at categories as generalisations of monoids. Let $\textbf{C}$ be a set equipped with a partial binary operation which we
shall denote by $\cdot$ or by concatenation. If $x, y\in \textbf{C}$ and the product $x \cdot y$
is defined we write $\exists x \cdot y$. An element $e \in \textbf{C}$ is called an \emph{identity} if
$\exists e \cdot x$ implies $e \cdot x = x$ and $\exists x \cdot e$ implies $x \cdot e = x$. The set of identities
of $\textbf{C}$ is denoted by $\textbf{C}_0$. The pair $(\textbf{C}, \cdot )$
is said to be a \emph{category} if the following axioms hold:\\

(C1): $x \cdot (y \cdot z)$ exists if, and only if, $(x \cdot y) \cdot z$ exists, in which case they are equal.

(C2): $x \cdot (y \cdot z)$ exists if, and only if, $x \cdot y$ and $y \cdot z$ exist.

(C3): For each $x \in \textbf{C}$ there exist identities $e$ and $f$ such that $\exists x\cdot e$ and $\exists f \cdot x$.\\
\\
It is convenient to write $xy$ instead of $x\cdot y$. From axiom (C3), it follows that the identities $e$ and $f$ are uniquely determined by $x$. We write $e = \r(x)$ and $f = \d(x)$.  We call $\d(x)$ the \emph{domain} of $x$ and $\r(x)$  the \emph{range}\footnote{Note that in \cite{hmark} Lawson used the converse notion, that is, $f = \r(x)$ and $e = \d(x)$ as he composed the functions from right to left.} of $x$. Observe that $\exists x · y$ if, and only if, $\r(x) = \d(y)$; in which case $\d(xy)=\d(x)$ and $\r(xy)=\r(y)$. The elements of $\textbf{C}$ are called \emph{homomorphisms}.\par
 \medskip
A  \emph{subcategory} $\mathbf{T}$ of a category $\mathbf{C}$ is a collection of some of the identities and some of the homomorphisms of  $\mathbf{C}$ which include with each homomorphism, $a$, both $\d(a)$ and $\r(a)$, and with each composable pair of homomorphisms in $\mathbf{T}$, their composite. 
In other words,  $\mathbf{T}$ is a category in its own right.\\
\\

The two definitions are equivalent. The first one  can  be easily turned into the second one and vice versa.  A homomorphism $a$ is said to be an \emph{isomorphism} if there exists an element $a^{-1}$ such that $\r(a)=a^{-1}a$ and $\d(a)=aa^{-1}$.\par
 \medskip
A \emph{groupoid} $\mathbb{G}$ is a category in which every element is an isomomorhism.  A group may be thought of as a one-object groupoid.  A category is  \emph{connected} if  for each pair of identities $e$ and $f$  there is a homomorphism  with domain $e$ and range $f$.  Connected groupoids are known as \emph{Brandt groupoids}.\par
 \medskip
If $\mathbb{G}$ and $\mathbb{P}$ are categories, then $\varphi:\mathbb{G} \longrightarrow \mathbb{P}$ is a \emph{homomorphis} if $\exists xy$ implies that $(xy)\varphi=(x\varphi) (y\varphi)$ and for all $x\in \mathbb{G}$ we have that $(\d(x))\varphi=\d(x\varphi)$ and $(\r(x))\varphi=\r(x\varphi)$. In case where $\mathbb{G}$ and $\mathbb{P}$ are groupoids we have that $\varphi:\mathbb{G} \longrightarrow \mathbb{P}$ is a homomorphism  if $\exists xy$ implies that $(xy)\varphi=(x\varphi) (y\varphi)$ and so $x^{-1}\varphi=(x\varphi)^{-1}$.\par
\medskip
The following lemma gives useful properties of groupoids which will be used without further mention. Proofs can be found in \cite{Ivan}. 

\begin{lemma}
Let $\mathbb{G}$ be a groupoid. Then for any $x,y \in \mathbb{G}$ we have \\

$(i)$ For all $x \in \mathbb{G}$  we have $\r(x^{-1})=\d(x)$ and $\d(x^{-1})=\r(x)$.

$(ii)$ If $\exists xy$, then $x^{-1}(xy)=y$ and $(xy)y^{-1}=x$ and $(xy)^{-1}=y^{-1}x^{-1}$.

$(iii)$ $(x^{-1})^{-1}=x$ for any $x\in \mathbb{G}$.
\end{lemma}

From now on we shall adopt the second definition of categories. In other words, we regard  categories as a generalisation of monoids.

 \begin{proposition}\label{prim2nd}\cite{Lawson98}
Let $G$ be a group and $I$ a non-empty set. Define a partial product on $I\times G\times I$ by $(i,g,j)(j,h,k)=(i,gh,k)$ and undefined in all other cases. Then $I\times G\times I$ is a connected groupoid, and every connected groupoid is isomorphic to one constructed in this way.
\end{proposition}

A \emph{Brandt semigroup} is a completely 0-simple inverse semigroup. By Theorem II.3.5 in \cite{Lawson98} every Brandt semigroup is isomorphic to  $B(G,I)$ for some group  $G$ and non-empty set  $I$ where  $B(G,I)$ is constructed as follows:\\ 
\\
As a set  $B(G,I)=(I \times G \times I)\cup \{0\},$ the binary operation is defined by

\begin{displaymath}
 (i,a,j)(k,b,l)= \left\{ \begin{array}{ll}
(i,ab,l), & \textrm{if \ $j= k$};\\
0, &  \textrm{else}
\end{array} \right.  
\end{displaymath} 
and
\[(i,a,j)0=0(i,a,j)=00=0.\]

 In \cite{clifford} it is shown that if we adjoint $0$ to  a Brandt groupoid $B$, defining $xy=0$ if $xy$ is undefined in $B$, we get a Brandt semigroup $B^0$.\\
 \\
 Let  $\{S_{i} : i \in I \}$ be  a family of disjoint semigroups with zero, and put  $S^*_i=S\setminus\{0\}$. Let  $ S=\bigcup_{i\in I} S^*_i\cup 0$  with the multiplication 
 \begin{displaymath}
a*b = \left\{ \begin{array}{ll}
ab, \ \mbox{if} \ a, b \in S_{i} \ \mbox{for some} \  i \ \mbox{and} \;  ab\neq 0 \ \mbox{in} \ S_{i};\\
\quad 0, \hspace{2cm} \mbox{else.}
\end{array} \right.\end{displaymath} 
 With this multiplication  $S$ is a semigroup called a \emph{0-direct union} of the $S_{i}$.\par
 \medskip
An  inverse semigroup   $S$ with zero  is a \emph{primitive inverse} semigroup if all its nonzero
idempotents are primitive, where an idempotent  $e$ of  $S$ is called \emph{primitive} if  $e\neq 0$ and  $f\leq e$ implies  $f=0$ or  $e=f$. 
Note that every Brandt semigroup is a primitive inverse semigroup.

\begin{theorem}\cite{Lawson98}
Brandt semigroups are precisely the connected groupoids with a zero adjoined, and every primitive inverse semigroup with zero is a 0-direct union of Brandt semigroups.
\end{theorem}

Notice that a groupoid is a disjoint union of its connected components.
\begin{theorem}\label{groupoidprimitive}\cite{Lawson98}
Let $\mathbb{G}$ be a groupoid. Suppose that $0\notin \mathbb{G}$ and put $\mathbb{G}^0=\mathbb{G}\cup \{0\}$. Define a binary operation on $\mathbb{G}^0$ as follows: if $x,y\in G$ and $\exists x.y$ in the groupoid $\mathbb{G}$, then $xy=x.y$; all other products in $\mathbb{G}^0$ are $0$. With this operation $\mathbb{G}^0$ is a primitive inverse semigroup.\end{theorem}

\begin{theorem}\cite{Lawson98}
Let $S$ be an inverse semigroup with zero. Then $S$ is primitive if, and only if, it is isomorphic to a groupoid with zero adjoined.\end{theorem}

 In \cite{clifford}, it is shown that every primitive inverse semigroup with zero is a 0-direct union of Brandt semigroups.\\
 \\

An \emph{ordered groupoid} $(\mathbb{G},\leq)$ is a groupoid $\mathbb{G}$ equipped with a partial order $\leq$ satisfies the following axioms:\\
 
(\textbf{OG1}) If $x\leq y$ then $x^{-1}\leq y^{-1}$.

(\textbf{OG2}) If $x\leq y$ and $x'\leq y'$ and the products $xx'$ and $yy'$ are defined then $xx'\leq yy'$.

(\textbf{OG3}) If $e\in \mathbb{G}_0$ is such that $e\leq \textbf{\d}(x)$ there exists a unique element $(x|e)\in \mathbb{G}$, called the \emph{restriction} of $x$ to $e$, such that $(x|e)\leq x$ and $\textbf{\d}(x|e)=e$.

(\textbf{OG3})$^*$ If $e\in \mathbb{G}_0$ is such that $e\leq \textbf{\r}(x)$ there exists a unique element $(e|x)\in \mathbb{G}$, called the \emph{corestriction} of $x$ to $e$, such that $(e|x)\leq x$ and $\textbf{\r}(e|x)=e$.\\
\\
In fact, it is shown in \cite{Lawson98} that axiom (\textbf{OG3})$^*$ is a consequence of the other axioms. \par
\medskip 
A  partially ordered set  $X$ is called a \textit{meet semilattice} if, for every $ x,y\in  X$, there is a greatest lower bound $x\wedge y$.  An ordered groupoid is \emph{inductive} if the partially ordered set of identities forms a meet-semilattice. An ordered groupoid $\mathbb{G}$ is said to be \emph{$\ast$-inductive} if each pair of identities that has a lower bound has a greatest lower bound. We can look at any inverse semigroup as an inductive groupoid; the order is the natural order and the multiplication is the usual multiplication.\par
\bigskip 
We shall now describe the relationship between inverse semigroups and  inductive groupoids. We begin with the following definition. 
\begin{definition}\normalfont
For an arbitrary inverse semigroup $S$, the \emph{restricted product} (also called the $\lq$trace product') of elements $x$ and $y$ of $S$ is $xy$ if $x^{-1}x=yy^{-1}$ and undefined otherwise.
\end{definition}

Let $S$ be an inverse semigroup with the natural partial order $\leq$. Define a partial operation $\circ$ on $S$ as follows:
\[x\circ y \ \mbox{defined iff}\ x^{-1}x=yy^{-1}\] in which case $x\circ y =xy$. Then $\mathbb{G}(S)=(S,\circ)$ is a groupoid and $\mathcal{G}(S)=(S,\circ, \leq)$ is an inductive groupoid with $(x|e)=xe$ and $(e|x)=ex$ and $e=x^{-1}xyy^{-1}$.

\begin{definition}\normalfont
Let $(\mathbb{G},.,\leq)$ be an ordered groupoid and let $x, y \in \mathbb{G}$ are such that $e = \textbf{\r}(x)\wedge \textbf{\d}(y)$ is defined. Then the \emph{pseudoproduct} of $x$ and $y$ is defined as follows:
\[x \otimes y = (x|e)(e|y).\] \end{definition}

If $(\mathbb{G},.,\leq)$ is an inductive groupoid, then $\mathcal{S}(\mathbb{G})=(\mathbb{G},\otimes)$ is an inverse semigroup having the same partial order as $\mathbb{G}$ such that the inverse of any element in $(\mathbb{G},.,\leq)$ coincides with the inverse of the same element in $\mathcal{S}(\mathbb{G})$. The pseudoproduct is everywhere defined in $\mathcal{S}(\mathbb{G})$ and coincides with the product $\cdot$ in $\mathbb{G}$ whenever $\cdot$ is defined, that is,  if $\exists x\cdot y$, then  $x \otimes y =x\cdot y$.\par
\medskip
 It is noted in \cite{Lawson98} that in an inductive groupoid $\mathbb{G}$, for $a\in \mathbb{G}$ and $e\in \mathbb{G}_0$ with $e\leq \r(a)$, the corestriction $e|a$ is given by $(e|a)=(a^{-1}|e)^{-1}$. By Linking this with the inverse semigroup which associated to $\mathbb{G}$, we present a short proof in the following lemma.

\begin{lemma}
Let $(\mathbb{G},.,\leq)$ be an inductive groupoid associated to an inverse semigroup $(\mathbb{G},\otimes)$. If $a\in \mathbb{G}$ and $e\in \mathbb{G}_0$ with $e\leq \r(a)$, then $(a^{-1}|e)^{-1}=(e|a)$. \end{lemma}
\begin{proof}
First we show that $(a^{-1}|e)$ exists. As $e\leq \r(a)$  and $\r(a)=\d(a^{-1})$ we have that $e\leq \d(a^{-1})$. Hence by (\textbf{OG3}) $(a^{-1}|e)$  exists. To show that $(a^{-1}|e)$ is the inverse of $(e|a)$. We note that
\[(a^{-1}|e) (e|a) (a^{-1}|e) =(a^{-1}e)(ea)(a^{-1}e)=a^{-1}e=(a^{-1}|e).\]
Also, 
\[(e|a)(a^{-1}|e)(e|a) =(ea)(a^{-1}e)(ea)=ea=(e|a).\]
 
\end{proof} 
  
We recall that a semigroup  $Q$ with zero is defined to be \emph{categorical at}  $0$  if whenever  $a,b,c \in Q$ are such that  $ab\neq 0$ and  $ bc \neq 0$, then  $abc\neq 0$. The set of non-zero elements of a semigroup $S$ will be denoted by $S^*$.\par
\medskip
Let $Q$ be an inverse semigroup which is categorical at zero. Define a partial binary operation $\circ$ on $Q^*$ by
\[a\circ b = \begin{cases}
ab, & \mbox{if} \ a^{-1}a=bb^{-1}; \\
\mbox{undefined}, & \mbox{otherwise}.  
\end{cases}\]
 It is easy to see that (C1) holds. Assume that $a\circ b$ and $b \circ c$ are defined in $Q^*$ so that $ab\neq 0$ and $bc\neq 0$ in $Q$. As $Q$ categorical at $0$ we have $abc\neq 0$ so that $a\circ(b\circ c)$ is defined in $Q^*$. On the other hand, if $a\circ(b\circ c)$ exists in $Q^*$, then $b^{-1}b=cc^{-1}$ and $a^{-1}a=(bc)(bc)^{-1}=bcc^{-1}b^{-1}=bb^{-1}$. Hence $a\circ b$ and $b\circ c$ exist. Thus (C2) holds. For any $a\in Q^*$ the identities $\d(a)=aa^{-1}$ and $\r(a)=a^{-1}a$  satisfy (C3). Hence $Q^*$ is a category and  any element $a$ in  $Q^*$ has the same inverse $a^{-1}$ as in $Q$. We have  
 
 \begin{lemma}\label{catatzerocat}
 Let $Q$ be an inverse semigroup with zero. If $Q$ categorical at $0$, then $Q^*=Q\setminus\{0\}$ is a groupoid.
 \end{lemma}

 Following \cite{Greenrela}, we define Green's relations on an ordered  groupoid $(\mathbb{G},\cdot,\leq)$. First we define  useful subsets of $\mathbb{G}$.\\
 \\
  For a subset $H$ of $\mathbb{G}$, we define $(H]$  as follows:

\[(H] = \{t \in \mathbb{G}: t \leq h \ \mbox{for some} \ h \in H\}.\]
For $a,b\in \mathbb{G}$, put \[\mathbb{G}a=\{xa:x\in \mathbb{G} \ \mbox{and}\ \exists xa\}.\]
 We define  $a\mathbb{G}$ and $a\mathbb{G}b$ similarly.\\
\\

A nonempty subset $I$ of $\mathbb{G}$ is called a \emph{right} (\emph{left}) \emph{ideal} of $\mathbb{G}$ if 

(1) $I\mathbb{G}\subseteq I$ ($\mathbb{G}I \subseteq  I$)

(2) if $a \in I$ and $ b \leq a$, then $b \in I$.\\
\\
We say that $I$ is  an \emph{ideal} of $\mathbb{G}$ if it is both a right and a left ideal of $\mathbb{G}$. We denote by $R(a), \ L(a), \ I(a)$
the right ideal, left ideal, ideal of $\mathbb{G}$, respectively, generated by $a$ ($a \in \mathbb{G}$). For each $a \in \mathbb{G}$, we have
\[R(a) = (a\cup a\mathbb{G}], \ L(a) = (a\cup \mathbb{G}a]\ \mbox{and} \ I(a) = (a\cup a\mathbb{G} \cup \mathbb{G}a\cup \mathbb{G}a\mathbb{G}].\]
For an ordered groupoid $\mathbb{G}$, the \emph{Green's relations} $\mathcal{R},\mathcal{L}$ and $\mathcal{J}$ defined on  $\mathbb{G}$ by
\[a\,\mathcal{R}\,b\ \Longleftrightarrow \ R(a) = R(b);\]
\[a\,\mathcal{L}\,b\ \Longleftrightarrow \ L(a) = L(b);\]
\[a\,\mathcal{J}\,b\ \Longleftrightarrow \ J(a) = J(b).\]

It is straightforward to show that $a\mathbb{G}=\d(a)\mathbb{G}$ ($\mathbb{G}a=\mathbb{G}\r(a)$) for all $a$ in $\mathbb{G}$. 
\begin{lemma}
Let $\mathbb{G}$ be an ordered groupoid and let $a,b \in \mathbb{G}$. Then

$(1)$ $a\,\mathcal{R}\,b \Longleftrightarrow \d(a) = \textbf{\d}(b)$.

$(2)$ $a\,\mathcal{L}\,b \Longleftrightarrow \textbf{\r}(a) = \textbf{\r}(b)$.
\end{lemma}
\begin{proof}
Suppose that $R(a) = R(b)$ it is clear that if $a=b$ we have that $a\,\mathcal{R}\,b$. If $a\neq b$, then $a\in R(b)$ so that $a\in (b\cup b\mathbb{G}] = \{t \in \mathbb{G}: t \leq h \ \mbox{for some} \ h \in b\cup b\mathbb{G}\}$. It is easy to see that $aa^{-1}\leq bb^{-1}$. Hence $\d(a)\leq \d(b)$. Similarly, we can show that $\d(b)\leq \d(a)$. Thus $\d(a)=\d(b)$.\par
Conversely, suppose that $\d(a)=\d(b)$. Let $x\in R(a)=(a\cup a\mathbb{G}]$ so that $x\leq h$ for some $h\in a\cup a\mathbb{G}$ so that $h=a$ or $h\in a\mathbb{G}=\d(a)\mathbb{G}=\d(b)\mathbb{G}=b\mathbb{G}$. In the latter case, it is clear that $x\in R(b)$. In the former case, $\d(h)=\d(a)$ and as $x\leq h$ we have that $xx^{-1}\leq hh^{-1}=\d(h)=\d(a)=\d(b)$ and so $x\leq \d(b)x$. Since $\d(b)x\in b\mathbb{G}\subseteq b\cup b\mathbb{G}$ we have that $x\in R(b)$ and so $R(a)\subseteq R(b)$. Similarly, $R(b)\subseteq R(a)$. Thus $R(b)= R(a)$ as required.
\end{proof}

%********\emph{semigroupoid} is a category without identity, also known as \emph{quiver}.******
%Green 's relations on a semigroupoid G are defined on $\Hom C$ as for semigroups. Moreover, in [PhD] was shown $\mathcal{R}\circ \mathcal{L}= \mathcal{L}\circ \mathcal{R}=\mathcal{D}$.

\section{Left orders in groupoids}\label{leftorder}

In this section consider the relationship between  left orders in an inductive groupoid  $\mathcal{G}$ and left I-orders in $\mathcal{S}(\mathcal{G})$.  We give a characterisation of left orders in groupoids. By using the second definition of categories we prove the category version of theorem due to Ore-Dubreil mentioned in the introduction.\\
\\

A category is said to be \emph{right} (\emph{left}) \emph{cancellative} if $\exists x\cdot a,\exists y\cdot a \ (\exists a\cdot x,\exists a\cdot y) $ and $xa=ya$ implies $x=y$ ($ax=ay$ implies $x=y$). A \emph{cancellative category} is one which is both left and right cancellative.\par
\medskip
Following \cite{hmark}, a category $\mathbf{C}$ is said to be \emph{right reversible}\footnote{ In \cite{hmark} Lawson used  $\textbf{\d}(a)=\textbf{\d}(b)$ instead of  $\textbf{\r}(a)=\textbf{\r}(b)$.} if for all $a,b\in \mathbf{C}$, with $\textbf{\r}(a)=\textbf{\r}(b)$, there exist $p,q\in \mathbf{C}$ such that $pa=qb$. In diagrammatic terms this just
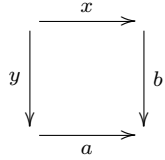
\[
\xymatrix @R=.5in @C=.5in{
&  \ar[d]_y \ar[r]^x
& \ar[d]^b \\
&  \ar[r]_a &  }
\] 

Let $\mathbf{C}$ be a category and $a,b\in \mathbf{C}$ such that $\d(a)=\d(b)$  we say that $a$ and $b$ have a \emph{pushout}, if $ax=by$ for some $x,y\in \mathbf{C}$. 
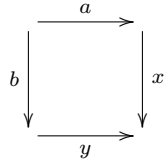
\[
\xymatrix @R=.5in @C=.5in{
&  \ar[d]_b \ar[r]^a
& \ar[d]^x \\
&  \ar[r]_y &  }
\] 
\begin{remark}\normalfont If a category $\textbf{C}$ is a left order in a groupoid $\mathbb{G}$, then any element in $\mathbb{G}$ has the form $a^{-1}b$.  It is clear that $\textbf{\d}(a)=\textbf{\d}(b)$ and $a^{-1}\,\mathcal{R}\,a^{-1}b\,\mathcal{L}\,b$. If any two elements   in $\textbf{C}$ have a  pushout, then $\textbf{C}$ is a right order in $\mathbb{G}$. Hence $\textbf{C}$ is an order in $\mathbb{G}$. We have the following diagram
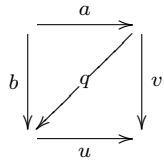
\[
 \xymatrix @R=.5in @C=.5in{
       \ar[r]^a \ar[d]_b  &   \ar[d]^v \ar[ld]|q &  \\
       \ar[r]_u &   &  \\
      & {} &  }\]\end{remark}

\begin{lemma}\label{fromindtoinv} A category $\textbf{C}$ is a left order in an inductive groupoid $\mathcal{G}$ if and only if $(\textbf{C},\otimes)$ is a left I-order in $(\mathcal{G},\otimes)$.\end{lemma}
\begin{proof}
Suppose that $\textbf{C}$ is a left order in $\mathcal{G}$. For any  $q\in \mathcal{G}$, there are $a,b\in \textbf{C}$ such that for $e=aa^{-1}bb^{-1}$ we have
\[\begin{array}{rcl}q&=&a^{-1}b  \\ &=&a^{-1}aa^{-1}bb^{-1}b\\ &=&(a^{-1}e)(eb)\\ &=&(ea)^{-1}(eb)\\ &=&(e|a)^{-1}(e|b)\\ &=& (a^{-1}|e)(e|b)\\ &=&a^{-1}\otimes b.\end{array}\]
It is clear that $(\textbf{C},\otimes)$ is a subsemigroup of $(\mathcal{G},\otimes)$. The converse follows by reversing the argument.
\end{proof}

\begin{lemma}\label{fromindtoinv}Let  $S$ be a semigroup which is a straight left I-order in an inverse semigroup $Q$. On the set $Q$ define a partial product $\circ$. Then $(S\cup E(Q),\circ)$ is a left order in $(Q,\circ)$.\end{lemma}
\begin{proof}
Suppose that $S$ is a straight left I-order in $Q$. For any  $q\in Q$, there are $c,d\in S$ such that $q=c^{-1}d$ with $c\,\mathcal{R}\,d$ so that $cc^{-1}=dd^{-1}$. Hence $q=c^{-1}d$ is defined in $(Q,\circ)$ and $c,d\in S\cup E(Q)$. It is easy to see that $E(Q)=\{a^{-1}a:a\in S\}$. Let $a^{-1}a\in E(Q)$ for some $a\in S$ and let $b\in S$ such that $ba^{-1}a$ is defined in $(Q,\circ)$ so that $b^{-1}b=a^{-1}a$. Hence $b=bb^{-1}b=ba^{-1}a\in S$. Similarly, if $a^{-1}ab$ is defined, then $a^{-1}a=bb^{-1}$ and so $b=bb^{-1}b=a^{-1}ab\in S$. Thus $(S\cup E(Q),\circ)$ is a left order in $(Q,\circ)$.
\end{proof}
The following corollary is clear.
\begin{corollary}\label{fromindtoinv2}Let  $S$ be a semigroup which is a straight left I-order in an inverse semigroup $Q$. On the set $Q$ define a partial product $\circ$. Then $(S,\circ)$ is a q-left order in $(Q,\circ)$.\end{corollary}

The following lemmas give a characterisation for categories which are  left orders in  groupoids. The proofs of such lemmas are quite straightforward and it can be deduced from \cite{Gaberiel} and \cite{Tob}, but we give it for completeness.

\begin{lemma}\label{rightreversiblecond}Let  $\textbf{C}$ be left order in a groupoid $\mathbb{G}$. Then

$(i)$ $\textbf{C}$  is cancellative;

$(ii)$ $\textbf{C}$ is  right reversible;

$(iii)$ any element in $\mathbb{G}_0$ has the form $a^{-1}a$ for some $a\in \textbf{C}$. Consequently, $\textbf{C}_0=\mathbb{G}_0$.
\end{lemma}
\begin{proof}
$(i)$ This is clear.\\
$(ii)$ Let $a,b\in \textbf{C}$ with $\r(a)=\r(b)$ so that $ab^{-1}$ is defined in $\mathbb{G}$. Since $\mathbb{G}$ is a category of left quotients of $\textbf{C}$, we have that $ab^{-1}=x^{-1}y$ where $x,y\in \textbf{C}$ and $\textbf{\d}(x)=\textbf{\d}(y)$. Then \[xa=xab^{-1}b=xx^{-1}yb=yb.\]
$(iii)$ Let $e$ be an identity in $\mathbb{G}_0$. As $\textbf{C}$ is a left order in $\mathbb{G}$ we have that  $e=a^{-1}b$ for some $a,b \in \textbf{C}$ so that $\d(a)=\d(b)$. Since $e$ is identity and $\d(a)=\d(b)$ we have  \[a=ae=aa^{-1}b=\d(b)b=b.\]
Hence $e=a^{-1}a=\r(a)\in \textbf{C}_0$ so that $\mathbb{G}_0\subseteq \textbf{C}_0$. Thus $\textbf{C}_0=\mathbb{G}_0$.
\end{proof}

\begin{lemma}\label{equality}
Suppose that $\mathbb{G}$ is a groupoid of left quotients of $\textbf{C}$. Then for all $a,b,c,d\in \textbf{C}$ the following are equivalent:

$(i)$ \ $a^{-1}b=c^{-1}d$;

$(ii)$ \ there exist $ x,y \in \textbf{C}$ such that  $xa=yc$ and $xb=yd$;

$(iii)$ \  $\textbf{\r}(a)=\textbf{\r}(c), \textbf{\r}(b)=\textbf{\r}(d)$ and for all $x,y \in S$ we have  $xa=yc \Longleftrightarrow xb=yd$.
                     
\end{lemma}
\begin{proof}
$(i)\Longrightarrow (ii)$ Suppose that $a^{-1}b=c^{-1}d$ for $a,b,c,d \in \textbf{C}$ so that $\r(a)=\r(c)$ and $\r(b)=\r(d)$. By Lemma~\ref{rightreversiblecond},  $\textbf{C}$ is right reversible and so there are elements $x,y\in \textbf{C}$ such that $xa=yc$. As, $\r(x)=\d(a)$ and  $\r(y)=\d(c)$ we have
\[ac^{-1}=x^{-1}xac^{-1}=x^{-1}ycc^{-1}=x^{-1}y.\]
Since $\d(a)=\d(b)$ and $\d(c)=\d(d)$ we have 
\[ca^{-1}=ca^{-1}bb^{-1}=cc^{-1}db^{-1}=db^{-1}.\]
Hence $db^{-1}=ca^{-1}=y^{-1}x$. As $\d(x)=\d(y)$ and $\r(b)=\r(d)$ we have that $xb=yd$.\\
\\
$(ii)\Longrightarrow (iii)$. It is clear that $\textbf{\r}(a)=\textbf{\r}(c)$ and $\textbf{\r}(b)=\textbf{\r}(d)$. Let $xa=yc$ and $xb=yd$. Suppose that $ta=rc$ for all $t,r\in \textbf{C}$. We have to show that $tb=rd$. By Lemma~\ref{rightreversiblecond},  $\textbf{C}$ is right reversible and cancellative.
Hence since $\r(y)=\r(r)$ and $\textbf{C}$, it follows that $ky=hr$ for some $k,h \in \textbf{C}$. Now,
\[kxa=kyc=hrc=hta,\]
cancelling in $\textbf{C}$ gives  $kx=ht$. Then 
\[htb=kxb=kyd=hrd,\]
again cancelling in $\textbf{C}$ gives $tb=rd$. Similarly, $tb=rd$ implies $ta=rc$, as required. \\
\\
$(iii)\Longrightarrow (i)$. Since $\textbf{C}$ is right reversible we have that $ta=rc$ for some $t,r\in \textbf{C}$  so that $tb=rd$. Then 
\[ac^{-1}=t^{-1}r=bd^{-1},\]
so that \[a^{-1}b=a^{-1}bd^{-1}d=a^{-1}ac^{-1}d=c^{-1}d,\]
as required.\end{proof}

 Lawson  has deduced  the following theorem from \cite{Gaberiel}. He has called the groupoid $\mathbb{G}$ in such a theorem a \emph{groupoid of fractions} of $\textbf{C}$.

\begin{theorem}\label{Markproof} \cite{hmark}
Let $\textbf{C}$ be a right reversible cancellative category. Then $\textbf{C}$ is a subcategory of a groupoid $\mathbb{G}$ such that the following three conditions hold:

$(i)$ \ $\textbf{C}_0=\mathbb{G}_0$.

$(ii)$ Every element of $\mathbb{G}$ is of the form $a^{-1}b$ where $a,b\in \textbf{C}$.

$(iii)$ $a^{-1}b=c^{-1}d$ if and only if there exist $x,y\in \textbf{C}$ such that $xa=yc$ and $xb=yd$.\end{theorem}

\begin{proof}
Our proof is basically the same as the proof given by Tobais \cite{Tob} in the case of categories as a collections of objects and homomorphisms, but our presentation is slightly different as we shall use the second definition of categories.\\
\\
From Lemmas ~\ref{rightreversiblecond} and ~\ref{equality},  $(i)$ and $(iii)$ are clear.\\
To prove $(ii)$  suppose that $\textbf{C}$ is right reversible and cancellative. We aim to construct a groupoid $\mathbb{G}$ in which $\textbf{C}$ is embedded as a left order in $\mathbb{G}$. This  construction is based on ideas by Tobias \cite{Tob} and Ghroda \cite{GG}. Let
\[\widetilde{\mathbb{G}}=\{(a,b)\in \textbf{C}\times \textbf{C}: \textbf{\d}(a)=\textbf{\d}(b)\}.\]

 Define a relation $(a,b)\sim (c,d)$ on $\widetilde{\mathbb{G}}$ by 
\[(a,b)\sim (c,d) \Longleftrightarrow \ \mbox{there exist}\ x,y \in \textbf{C}\ \mbox{such that} \ xa=yc\ \mbox{and}\ xb=yd.\]
We can  represent this relation by the following diagram
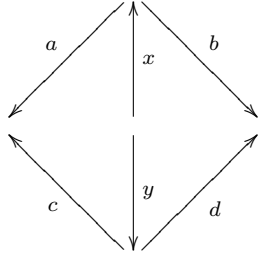
\[
\begin{xy}
\xymatrix @R=.6in @C=.6in{
&  \ar[dl]_a \ar[rd]^b & \\
  & \ar[u]_x \ar[d]^y
& \\  &\ar[lu]^c \ar[ru]_d &}
\end{xy}
\]
Notice that if $(a,b)\sim (c,d)$, then $\r(a)=\r(c)$ and $\r(b)=\r(d)$.
\begin{lemma}
The relation $\sim$ defined above is an equivalence relation.
\end{lemma}
\begin{proof}
It is clear that $\sim$ is symmetric and reflexive.  Let 
\begin{equation*}
(a,b) \sim (c,d)\sim (p,q),
\end{equation*}
where $(a,b),(c,d)$ and $(p,q)$ in $\widetilde{\mathbb{G}}$. Hence there exist $x , y ,\bar{x} , \bar{y}\in \textbf{\textbf{C}} $ such that  
\begin{equation*}
xa=yc  , \ xb=yd  \ \mbox{and} \ \bar{x}c=\bar{y}p, \ \bar{x}d=\bar{y}q.
\end{equation*}
To show that $\sim$ is transitive, we have to show that there are elements $z , \bar{z} \in \textbf{C}$ such that $za=\bar{z}p$ and $zb=\bar{z}q$. \newline
\newline
Since $\textbf{C}$ is right reversible and $\textbf{\r}(y)=\textbf{\r}(\bar{x})$ there are elements $s,t\in \textbf{C}$ such that $sy=t\bar{x}$. Hence \[sxa=syc=t\bar{x}c=t\bar{y}p.\]
Similarly, $sxb=t\bar{y}q$ as required.
\end{proof}

Let $[a,b]$ denote the $\sim$-equivalence class of $(a,b)$. On $\mathbb{G}=\widetilde{\mathbb{G}}/$$\sim$ we define a product as follows. Let $[a,b],[c,d]\in \mathbb{G}$. Their product is defined iff $\r(b)=\r(c)$. Define

\[[a,b][c,d] = \begin{cases}
[xa,yd], & \ \mbox{if}\  xb=yc \  \mbox{for some}\ x,y \in \textbf{C}; \\
\mbox{undefined}, & \mbox{otherwise},  
\end{cases}\] 
and so we have the following diagram
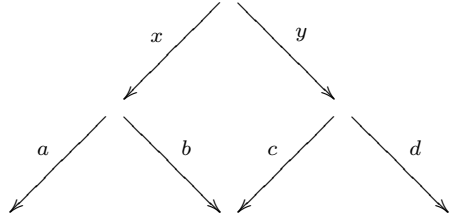
\[\begin{xy}
\xymatrix @R=.5in @C=.5in{
&&  \ar[dr]^y \ar[dl]_x & \\
& \ar[dl]_a\ar[dr]^b & & \ar[dl]_c\ar[dr]^d&&
 \\  & &  & & }
\end{xy}
\]
\begin{lemma}\label{welldefinedgroupoid} 
The multiplication is well-defined.
\end{lemma}
\begin{proof}
Suppose that $[a_{1},b_{1}]=[a_{2},b_{2}]$ and $[c_{1},d_{1}]=[c_{2},d_{2}]$ are
in $\mathbb{G}$. Then there are elements $x_{1} , x_{2} , y_{1} , y_{2}$ in $\textbf{C}$
such that 
\begin{equation*}
\begin{array}{rcl}
x_{1}a_{1}=x_{2}a_{2}, &  &  \\ 
x_{1}b_{1}=x_{2}b_{2}, &  &  \\ 
y_{1}c_{1}=y_{2}c_{2}, &  &  \\ 
y_{1}d_{1}=y_{2}d_{2}. &  & 
\end{array}%
\end{equation*}
Now, 
\begin{equation*}
[a_1,b_1][c_{1},d_1] = [wa_1,\bar{w}d_1] \; \mbox{and}
\;wb_{1}=\bar{w}c_{1} \end{equation*}
for some $w,\bar{w}\in \textbf{C}$ and 
\begin{equation*}
[a_{2},b_{2}][c_{2},d_{2}] = [za_2,\bar{z}d_2]  \; \mbox{and}
\; zb_{2}=\bar{z}c_{2}\end{equation*}
for some $z,\bar{z}\in \textbf{C}$.\\
It is easy to see that $[a_1,b_1][c_{1},d_1]$ is defined if and only if $[a_{2},b_{2}][c_{2},d_{2}]$ is defined.
 We have to prove that $[wa_{1},\bar{w}d_{1}]=[za_{2},\bar{z}d_{2}]$, that is, 
\begin{equation*}
xwa_{1}=yza_{2}\  \mbox{and} \ x\bar{w}d_{1}=y\bar{z}d_{2}, \; \mbox{for some} \; x , y \in \textbf{C}. 
\end{equation*}
Since $wb_{1}$ is defined and $\d(a_1)=\d(b_1)$ we have that $wa_1$ is defined and $\r(a_1)=\r(wa_1)$. Similarly, $za_2$ is defined and $\r(a_2)=\r(za_2)$. Hence $\r(wa_1)=\r(wa_2)$, by the right reversibility of $\textbf{C}$ there are elements $x,y \in \textbf{C}$ with $xwa_1=yza_2$. It remains to show that $x\bar{w}d_{1}=y\bar{z}d_{2}$. By Lemma~\ref{equality}, $xwb_1=yzb_2$ and as $wb_{1}=\bar{w}c_{1}$ and $zb_{2}=\bar{z}c_{2}$ we have that $x\bar{w}c_{1}=y\bar{z}c_{2}$ and so $x\bar{w}d_{1}=y\bar{z}d_{2}$, again by Lemma~\ref{equality}. 
\end{proof}

\begin{lemma}\label{associative}
The multiplication is associative.
\end{lemma}
\begin{proof}
Let $[a,b],[c,d],[p,q]\in \mathbb{G}$ and set 
\begin{equation*}
X=([a,b][c,d])[p,q] = [xa,yd][p,q]   
\end{equation*}
where $xb=yc$ for some $x,y\in \textbf{C}$ and 
\begin{equation*}
Y=[a,b]([c,d][p,q]) = [a,b][\bar{x}c,\bar{y}q] 
\end{equation*}
where $\bar{x}d=\bar{y}p$ for some $\bar{x},\bar{y}\in \textbf{C}$. It is clear that $X$ is defined if and only if $Y$ is defined.
 We assume that $[a,b][c,d]$ and  $[c,d][p,q]$ are defined. Then for some $x,y,\bar{x},\bar{y} \in \textbf{C}$ we have
\[\begin{array}{rcl}
X&=&[xa,yd][p,q]\\ &=&[sxa,rq]  \end{array}\]
where $syd=rp$ for some $ s,r \in \textbf{C}$.
 
\[\begin{array}{rcl}
Y&=&[a,b][\bar{x}c,\bar{y}q]\\ &=&[\bar{s}a,\bar{r}\bar{y}q]  \end{array}\]
where $\bar{s}b=\bar{r}\bar{x}c$ for some $\bar{s},\bar{r} \in \textbf{C}$.
We have to show that \[X=[sxa,rq]=[\bar{s}a,\bar{r}\bar{y}q]=Y.\] Then by definition we need to show that
\begin{equation*}
wsxa=\bar{w}\bar{s}a \ \mbox{and} \ wrq=\bar{w}\bar{r}\bar{y}q
\end{equation*}
for some $w , \bar{w} \in \textbf{C}$. By cancellativity in $\textbf{C}$ this equivalent to $wsx=\bar{w}\bar{s}$ and $wr=\bar{w}\bar{r}\bar{y}$. \par
\medskip
Since $xb$ and $sx$ are defined we have that $sxb$ is defined and as $\bar{s}b$ is defined so that $\r(\bar{s}b)=\r(sxb)$. By right reversibility of $\textbf{C}$ we have that $wsxb=\bar{w}\bar{s}b$  for some $w,\bar{w} \in \textbf{C}$. By  cancellativity in  $\textbf{C}$ we get $wsx=\bar{w}\bar{s}$. \newline
\newline
Now, since $wsxb=\bar{w}\bar{s}b , \ \bar{s}b=\bar{r}\bar{x}c$ and $xb=yc$ we have that $wsyc=\bar{w}\bar{r}\bar{x}c$. As $\textbf{C}$ is cancellative we have that $wsy=\bar{w}\bar{r}\bar{x}$ so that $wsyd=\bar{w}\bar{r}\bar{x}d$, but $syd=rp$ and $\bar{x}d=\bar{y}p$ so that $wrp=\bar{w}\bar{r}\bar{y}p$. Thus $wr=\bar{w}\bar{r}\bar{y}$ as required.
\end{proof}

For $[a,b] \in \mathbb{G}$ where $xa$ is defined in $\textbf{C}$ for some $x\in \textbf{C}$, it is clear that $[xa,xb]\in \mathbb{G}$ and $\d(x)xa=xa$ and $\d(x)xb=xb$. Hence we have

\begin{lemma}\label{supportive2} If $[a,b], [xa,xb] \in \mathbb{G}$, then $[xa,xb]=[a,b]$ for all $x\in \textbf{C}$ such that  $xa$ is defined in $\textbf{C}$.\end{lemma}
 \begin{lemma}\label{identities}
 The identities of $\mathbb{G}$ have the form $[a,a]$ where $a\in \textbf{C}$.\end{lemma}
 \begin{proof}
 Suppose that $e=[a,b]$ is an identity in $\mathbb{G}$ where $a,b\in \textbf{C}$. Let $[m,n]\in \mathbb{G}$ such that 
$[m,n][a,b]$ is defined and \[[m,n][a,b]=[m,n].\]
 Then $[xm,yb]=[m,n]$ for some $x,y \in \textbf{C}$ with $xn=ya$. Hence \[uxm=vm \  \mbox{and} \ uyb=vn\]  for some $u,v \in \textbf{C}$;  cancelling in $\textbf{C}$ gives that $ux=v$ so that $uyb=vn=uxn$. Again, by cancellativity, it follows that $xn=yb$ and as $xn=ya$ we have that $yb=xn=ya$. Using cancellativity  in $\textbf{C}$ once more we obtain $a=b$. Thus $e=[a,a]$. Similarly, if $[a,b][m,n]=[m,n]$ we have that $a=b$. It remains to show that the identity is unique. Suppose that \[[a,b][c,c]=[a,b][d,d]=[a,b]\] for some identities $[c,c],[d,d] \in \mathbb{G}$. Then by definition $[xa,yc]=[x'a,y'd]$ where $xb=yc$ and $x'b=y'd$ for some $x,y,x',y'\in \textbf{C}$. Hence $uxa=vx'a$ and $uyc=vy'd$. By definition of $\sim$ and Lemma~\ref{supportive2}, \[[c,c]=[yc,yc]=[y'd,y'd]=[d,d].\]
 As required. Similarly, $[c,c][a,b]=[d,d][a,a]=[a,b]$ implies that $[c,c]=[d,d]$.
 \end{proof}

 Suppose that $[a,b]\in \mathbb{G}$. Then as $\d(a)a=\d(a)a$ we have
 \[[a,a][a,b]=[\d(a)a,\d(a)b]=[a,b].\]
Similarly, $\d(b)b=\d(b)b$ whence
\[[a,b][b,b]=[\d(a)a,\d(a)b]=[a,b].\]
Hence $\d([a,b])=[a,a]$ and $\r([a,b])=[b,b]$.\\
 By the above argument and Lemma~\ref{associative}, the following lemma is clear.
 
 \begin{lemma} $\mathbb{G}$ is a category.\end{lemma}

If $[a,b]\in \mathbb{G}$, then it is clear that $[b,a] \in \mathbb{G}$, as $\d(a)b=\d(a)b$ we have  
 \[[a,b][b,a]=[\d(a)a,\d(a)a]=[a,a]=\d([a,b])\]
 Similarly, $[b,a][a,b]=[b,b]=\r([a,b])$. That is, $[b,a]$ is the inverse of $[a,b]$ in $\mathbb{G}$. Thus we have

  \begin{lemma} $\mathbb{G}$ is a groupoid.\end{lemma}
  
  Although from $(i)$ we know that $\textbf{C}_0=\mathbb{G}_0$, but in the following lemma we give a new proof depends on the structure of $\mathbb{G}$.
  \begin{lemma}  $\textbf{C}_0=\mathbb{G}_0$. \end{lemma}
 \begin{proof}
 Suppose that $e=[a,b]$ is an identity in $\mathbb{G}$ where $a,b\in \textbf{C}$. Then
 \[[a,b]=[\d(x),x]^{-1}[\d(y),y].\]
 As $[\d(x),x][\d(x),x]^{-1}$ is defined, so that $[\d(x),x][a,b]$ is defined and since $[a,b]$ is an identity we have that  $[\d(x),x][a,b]=[\d(x),x]$. Hence
 \[[\d(x),x]=[\d(x),x][a,b]=[\d(x),x][\d(x),x]^{-1}[\d(y),y]=[\d(y),y].\]
 Thus $x=y$ and so $ a=b$. But $[a,a]=[\r(a),\r(a)]=[\d(\r(a)),\r(a)]$ so that the identities in $\textbf{C}$ have the form $[\d(\r(a)),\r(a)]$ for any $a\in \textbf{C}$. 
 \end{proof}

 \begin{lemma}\label{embedding} The mapping $\theta: \textbf{C} \longrightarrow \mathbb{G}$ defined by $a\theta=[\d(a),a]$  is an embedding of $\textbf{C}$ in $\mathbb{G}$.
 \end{lemma}
\begin{proof}
It is clear that $\theta$ is well-defined. To show that $\theta$ is one-to-one, let $[\d(a),a]=[\d(b),b]$ so that $ua=vb$ and $u\d(a)=v\d(b)$ for some $u,v\in \textbf{C}$. Hence $a=b$.\\
\\
Let $a,b\in \textbf{C}$ such that $ab$ is defined. We have 
\[\begin{array}{rcl}a\theta b\theta &=& [\d(a),a][\d(b),b]\\ &=&[u\d(a),vb] \hspace{1.8cm} \mbox{where}\ ua=v\d(b) \ \mbox{for some} \ u,v \in \textbf{C} \\ &=&[u\d(a),uab] \hspace{1.6cm}  \mbox{as} \ ua=v\d(b)=v \\ &=&[\d(a),ab] \hspace{2.1cm} \mbox{by lemma~\ref{supportive2}}\\ &=&[\d(ab),ab] \hspace{1.9cm}  \mbox{as}\ \d(a)=\d(ab)\\  &=&(ab)\theta.\end{array}\]
Thus $\theta$ is a homomorphism.
\end{proof}

By Lemma~\ref{embedding}, we can consider  $\textbf{C}$ as a subcategory of $\mathbb{G}$. Let $[a,b]\in \mathbb{G}$  and  $a\theta=[\d(a),a]$, $b\theta=[\d(b),b]$. As $\d(a)=\d(b)$ we have that $a\theta=[\d(a),a]$ and $b\theta=[\d(a),b]$. Hence
 \[\begin{array}{rcl}(a\theta)^{-1} (b\theta) &=& [\d(a),a]^{-1}[\d(a),b]\\ &=& [a,\d(a)][\d(a),b]\\ &=& [ua,vb]\hspace{2cm}  \mbox{where} \ u=u\d(a)=v\d(a) =v\\ &=&[ua,ub] \hspace{2cm}  \mbox{by Lemma~\ref{supportive2}}\\ &=&[a,b] .\end{array}\]
Hence $\textbf{C}$ is a left order in  $\mathbb{G}$. This completes the proof of Theorem~\ref{Markproof}.
\end{proof}

 \begin{corollary}\label{catequ}
 A subcategory $\textbf{C}$ is  a left order in a groupoid $\mathbb{G}$ \textit{if and only if }   $\textbf{C}$ is right reversible and cancellative.
\end{corollary}
\begin{proof}
 If $\textbf{C}$ is a left order a groupoid $\mathbb{G}$, then by Lemmas~\ref{rightreversiblecond},  $\textbf{C}$ is right reversible and cancellative. Conversely, if $\textbf{C}$ is right reversible and cancellative, then by $(ii)$ in Theorem~\ref{Markproof}, $\textbf{C}$ is a left order a groupoid.
\end{proof}

\section{Uniqueness}\label{uniq} 
In this section we show that  a category $\textbf{C}$  has, up to isomomorphism,  at most one  groupoid of left I-quotients.

  \begin{theorem}
 Let $\textbf{C}$ be a left order in groupoid $\mathbb{G}$. If $\varphi$ is an embedding of  $\textbf{C}$ to a groupoid $\mathbb{T}$, then there is a unique embedding $\psi: \mathbb{G}\longrightarrow \mathbb{T}$ such that $\psi|_S=\varphi$.
 \end{theorem}
 
\begin{proof}
Define $\psi: \mathbb{G}\longrightarrow \mathbb{T}$ by 
\[(a^{-1}b)\psi=(a\varphi)^{-1}(b\varphi)\]
for $a,b,c,d \in \textbf{C}$. Suppose that \[a^{-1}b=c^{-1}d\] so that $xa=yc$ and $xb=yd$ for some $x,y \in \textbf{C}$, by Lemma~\ref{equality}. Hence \[x\varphi a\varphi=y\varphi c\varphi\ \mbox{and} \ x\varphi b\varphi=y\varphi d\varphi\] in $\textbf{C}\varphi$. Thus
 \[a\varphi c\varphi^{-1} =x\varphi^{-1} y\varphi= b\varphi d\varphi^{-1}\] so that \[a\varphi^{-1}b\varphi = c\varphi^{-1} d\varphi.\] It follows that $\psi$ is well-defined and 1-1. It remains for us to show that $\psi$ is a homomorphism. Let $a^{-1}b,c^{-1}d\in \mathbb{G}$ where $a,b,c,d \in \textbf{C}$. Now,
  \[\begin{array}{rcl}(a^{-1}bc^{-1}d)\psi&=&((xa)^{-1}(yd))\psi \\ &=&(xa)\varphi^{-1}(yd)\varphi\\ &=& a\varphi^{-1} x\varphi^{-1}y\varphi d\varphi,\end{array} \]
where $xb=yc$ for some $x,y \in \textbf{C}$. We have that $x\varphi b\varphi=y\varphi c\varphi$ and so $ b\varphi c\varphi^{-1}=x\varphi^{-1} y\varphi$. Hence
  \[\begin{array}{rcl}(a^{-1}bc^{-1}d)\psi &=& a\varphi^{-1} x\varphi^{-1}y\varphi d\varphi \\ &=& a\varphi^{-1} b\varphi c\varphi^{-1} d\varphi \\ &=&(a^{-1}b)\psi(c^{-1}d)\psi.\end{array} \]
 
 Finally, to see that $\psi$ is unique, suppose that $\theta : \mathbb{G}\longrightarrow \mathbb{T}$ is an embedding with $\theta|_S=\varphi$. Then for an element $a^{-1}b$ of $\mathbb{G}$, we have
\[(a^{-1}b)\theta = (a^{-1}\theta) (b\theta) = (a\theta)^{-1} (b\theta) = (a\varphi)^{-1} (\varphi) = (a^{-1}b)\psi \]
so that $\theta =\psi$.
 
 \end{proof}
 
The following corollary is straightforward.
 \begin{corollary}\label{coruniqu}
If a category $\textbf{C}$ is a left order in groupoids $\mathbb{G}$ and  $\mathbb{P}$,
then $\mathbb{G}$ and  $\mathbb{P}$ are isomorphic by an isomorphism which restricts to the identity map on $\textbf{C}$.
 \end{corollary}

\section{Primitive inverse semigroups of left I-quotients and groupoids of left quotients}\label{primitiveord}

In this section we are concerned with the relationship between primitive inverse semigroups of left I-quotients and groupoids of left quotients.\par
\medskip
 Let  $Q$ be a primitive inverse semigroup. By using the restriction product of elements of $Q$ we can associate a groupoid  $\mathbb{G}$ to $Q$. On the other hand, if $\mathbb{G}$ is a groupoid, then $\mathbb{G}^0$ is a primitive inverse semigroup.  In fact,  any primitive inverse semigroup is isomorphic to one constructed in this way. In particular, primitive inverse semigroups of left I-quotients and groupoids of left quotients are equivalent structures in the sense that each can be reconstructed from the other.  \par
\medskip
It is well-known that any primitive inverse semigroup is categorical at $0$.  We say that,  $S$ is \emph{0-cancellative} if  $b=c$ follows from  $ab=ac\neq 0$  and from  $ba=ca\neq0$. \par
\medskip
In \cite{GG} we studied left I-orders in  primitive inverse semigroups by using the relation  $\lambda$ on any semigroup with zero which defined as follows:
\[a\,\lambda\,b \ \mbox{if and only if} \ a=b=0 \ \mbox{or}\ Sa\cap Sb\neq 0. \]

The following theorem gives necessary and sufficient conditions for a semigroup to have a primitive inverse semigroup of left I-quotients.
 \begin{theorem}\label{primitiveq}\cite{GG}
 A semigroup $S$ is a left I-order in a primitive inverse semigroup $Q$ \textit{if and only if } $S$ satisfies the following conditions: 
 
$(A)$ \ $S$ is categorical at $0$; 

$(B)$ \ $S$ is 0-cancellative;

$(C)$ \ $\lambda$ is transitive; 
 
$(D)$ \ $Sa\neq 0$ for all $a\in S^*$.
\end{theorem}

Suppose that a semigroup $S$ has a primitive inverse semigroup of left I-quotients $Q$. By Proposition 2.4 of \cite{GG}, $S$   contains a zero. Define a partial binary operation $\circ$ on $\mathbb{G}=Q^*$ by
\[a\circ b = \begin{cases}
ab, & \mbox{if} \ ab\neq 0; \\
\mbox{undefined}, & \mbox{otherwise}.  
\end{cases}\] 
Notice that in $Q$, for any non-zero elements $a$ and $b$ in $Q$  we have that $ab\neq 0$ if and only if $a^{-1}a=bb^{-1}$. We have a groupoid $\mathbb{G}$ with $\mathbb{G}_0=E(Q)$ ordered by the natural partial order on $Q$. Hence $\mathbb{G}$  is inductive.\par
 \medskip
If a semigroup $S$ is a left I-order in  $Q$, then $S$ satisfies the conditions of Theorem~\ref{primitiveq}. Since  $Q$ categorical at $0$, it follows that $Q^*$ with the partial product is a groupoid, by Lemma~\ref{catatzerocat}. It is clear that any element in $Q^*$ can be written as $a^{-1}b$ where $a,b\in S^*$. Unfortunately, $S^*$ is not a subcategory of $Q^*$ as we do not insist on $S$ being full in $Q$. By adding $E(Q)$ to $S^*$ we have that $\textbf{C}=S^*\cup E(Q)$ is a subcategory of $Q^*$. For, let $a\in S^*$ and $f\in E(Q)$, if $af$ is defined in $ \textbf{C}$, then $af\neq 0$ in $Q$. By Lemma 2.1 in \cite{Gjhon} we have $a^{-1}a=f$ and so $a=aa^{-1}a=af$ so that $a=af\in \textbf{C}$. Similarly, if $ea$ is defined in $\textbf{C}$, then $a=ea$. Thus $\textbf{C}$ is a left order in $Q^*$. We have

 \begin{lemma}\label{fromprimitivetog}
If a semigroup $S$ is a left I-order in a primitive inverse semigroup $Q$, then $\textbf{C}=S^*\cup E(Q)$ is a left order in the groupoid $Q^*$.
 \end{lemma}
  
  \begin{corollary}\label{fromprimitivetoC}
If a semigroup  $S$ is a left I-order in a primitive inverse semigroup $Q$, then $\textbf{C}=S^*\cup E(Q)$ is a right reversible cancellative subcategory of the groupoid $\mathbb{G}=Q^*$ associated to $Q$. 
 \end{corollary}

 Now, in case we did not drop the zero from $Q$. Under the same multiplication we can show that $\textbf{C}=(S\cup E(Q),\circ)$ is a left order in $\mathbb{G}=(Q,\circ)$ where $\mathbb{G}$ is inductive and $\{0\}$ is an isolated  identity of $\mathbb{G}$, in the sense that it is not a domain or codomain of any hommorphism. Also, $\mathbb{G}^*=\mathbb{G}\setminus\{0\}$ is an *-inductive groupoid with $\mathbb{G}^*_0=E(Q)^*$.\par
 \medskip
 We have shown how to construct a groupoid of left quotients from a primitive inverse semigroup of left I-quotients. The rest of this section is devoted to showing how to construct  a primitive inverse semigroup of left I-quotients from a groupoid of left quotients.\par
 
 \bigskip
 Suppose that $\textbf{C}$ is a left order in a groupoid $\mathbb{G}$.   Define multiplication on $Q$ as:
 \[a b = \begin{cases}
a\circ b, & \mbox{if} \ \exists a\circ b \: \mbox{in} \ \mathbb{G};\\
0, & \mbox{otherwise}.  
\end{cases}\] 
 
By Theorem~\ref{catequ}, $\textbf{C}$   is right reversible and cancellative. Put $S=\textbf{C}\cup \{0\}$ and $Q=\mathbb{G}\cup \{0\}$. By Theorem~\ref{groupoidprimitive}, $Q$ is a primitive inverse semigroup. It is clear that $S$ is a subsemigroup of $Q$. Moreover, it is a left I-order in $Q$.  Thus we have
 
 \begin{lemma}\label{fromgroupidtop}
If a category $\textbf{C}$ is a left order in a groupoid $\mathbb{G}$, then $S=\textbf{C}\cup \{0\}$ is a left I-order in $Q=\mathbb{G} \cup\{0\}$.
 \end{lemma}

 Notice that $S$ in the above lemma is full, as $\textbf{C}$ and $\mathbb{G}$ have the same set of the identites.

\section{Connected groupoids of left quotients}\label{connectedored}

In this section we give necessary and sufficient conditions for a category to have a connected groupoid of left quotients. That is, we  specialise our result in the previous section to left orders in connected groupoids. In general a groupoid it might contains a zero, but a non-trivial connected groupoid can not have a zero. \par
 \medskip
 
Let $\textbf{C}$ be a category. We say that $\textbf{C}$ satisfies the \textit{connected condition} if,  for any $a,b\in \textbf{C}$ there exist $c,d\in \textbf{C}$ with $\d(c)=\d(d)$ such that  $\r(c)=\r(a)$  and  $\r(d)=\r(b)$. We regard the connected condition in diagrammatic terms as follows.
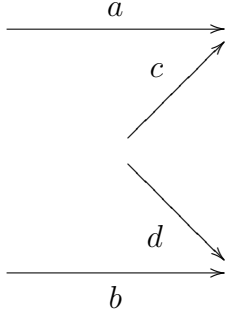
\[\Large\begin{xy}
 \xymatrix @R=.5in @C=.5in{ 
       \ar[rr]^a  &    &  \\
      & \ar[dr]_d \ar[ru]^c &  \\
     \ar[rr]_b   & {} &  }\end{xy}\]
     
Let  $\textbf{C}$ be  a left order in a groupoid $\mathbb{G}$. It is clear that if $\textbf{C}$ is connected, then $\mathbb{G}$ is connected. Also, if $\textbf{C}$ is connected, then it has the connected condition, but the converse is not true.
\begin{lemma}\label{connected}
A  right reversible, category $\textbf{C}$ is a left order in a  connected groupoid $\mathbb{G}$  if and only if $\textbf{C}$ satisfies the connected condition.\end{lemma}
\begin{proof}
Suppose that $\textbf{C}$ is a left order in a  connected groupoid $\mathbb{G}$. To show that $\textbf{C}$ satisfies the connected condition.  Suppose that $a^{-1}a$ and $b^{-1}b$ are two identity elements of $\mathbb{G}$ for some $a,b \in\textbf{C}$. Hence there is an isomomorphism $h$ in $\mathbb{G}$ such that $\d(h)=a^{-1}a$ and $\r(h)= b^{-1}b$. As $\textbf{C}$  is a left order in $\mathbb{G}$, we have that $h=s^{-1}t$ for some $s,t \in \textbf{C}$. Hence $\r(s)=\d(h)=\r(a)$ and $\r(t)=\r(h)=\r(b)$. Thus  the connected condition holds on $\textbf{C}$. \par
Conversely, suppose that $\textbf{C}$  satisfies the connected condition. By Corollary~\ref{catequ}, $\textbf{C}$ is a left order in a groupoid $\mathbb{G}$. We proceed to show that $\mathbb{G}$ is connected. Suppose that $e$ and $f$ are two identity elements of $\mathbb{G}$ so that $e=a^{-1}a$ and $f=b^{-1}b$ for some $a,b \in \textbf{C}$, by Lemma~\ref{rightreversiblecond}. By assumption there exist $c,d\in \textbf{C}$ with $\d(c)=\d(d)$ such that  $\r(b)=\r(d)$ and $\r(c)=\r(a)$. As  $\d(c)=\d(d)$ we have $c^{-1}d$ is defined in $\mathbb{G}$ and $\d(c^{-1}d)= c^{-1}c=a^{-1}a$ and  $\r(c^{-1}d)=d^{-1}d=b^{-1}b$. Thus $\mathbb{G}$ is connected.
\end{proof}  

%\begin{corollary}$\mathbb{G}$  is a connected groupoid if and only if $\mathbb{G}^0$  is a 0-bisimple inverse semigroup.??\end{corollary}

\begin{lemma} 
Let  $\mathbb{G}=\bigcup_{i\in I} \mathbb{G}_i$  be a groupoid where  $\mathbb{G}_i$  is a connected groupoid. If  $\textbf{C}$ is a left order in  $\mathbb{G}$, then  $\textbf{C}$  is a disjoint union of categories that are left orders in the connected groupoids $\mathbb{G}_i$'s.\end{lemma}
\begin{proof}
Suppose that $\textbf{C}$ is a left order in $\mathbb{G}$. Then every element of $\mathbb{G}$ can be written as $a^{-1}b$ where $a,b \in \textbf{C}$. As $\mathbb{G}$ is a disjoint union of $\mathbb{G}_i$'s we have that $a^{-1},b \in \mathbb{G}_i$ for some $i\in I$. Since $\mathbb{G}_i$ is a groupoid we have that $a,b \in \mathbb{G}_i$. Hence $\textbf{C}_i= \textbf{C}\cap \mathbb{G}_i\neq \emptyset$. It is clear that $\textbf{C}_i$ is a subcategory of $\mathbb{G}_i$. We claim that $\textbf{C}_i$ is a left order in $\mathbb{G}_i$. Let $q=a^{-1}b\in \mathbb{G}_i$ so that $a,b \in \textbf{C}\cap \mathbb{G}_i=\textbf{C}_i$. It is clear that $\textbf{C}$ is a disjoint union of the categories $\textbf{C}_i, i\in I$.  

\end{proof}
 
 It is well-known that Brandt semigroups are precisely the connected groupoids with a zero adjoined. On the other hand, from any Brandt semigroup we can recover a connected groupoid. We will be using the same technique that we used in the previous section to determine the relationship between  left I-orders in Brandt semigroups and left orders in connected groupoids.  Brandt semigroups are  primitive so that the most part of the task has been done in the previous section. 
 
\begin{theorem}\cite{GG}
 A semigroup $S$ is a left I-order in a Brandt semigroup $Q$ \textit{if and only if } $S$ satisfies the following
conditions: 

$(A)$ \ $S$ is categorical at $0$; 

$(B)$ \ $S$ is 0-cancellative; 
 
$(C)$ \ $\lambda$ is transitive; 

$(D)$ \ $Sa\neq 0$ for all $a\in S^*$;

$(E)$ \ for all $a , b \in S^*$ there exist $c , d \in S$ such that  $ca\,\mathcal{R^{*}}\,d\,\lambda\,b$.
\end{theorem}

In order to generalise the above theorem to the category version, we need to obtain a corresponding condition to (E) to make $\mathbb{G}$ connected, that is, we need to show that there is  an isomorphism between any two identities. \\
\\

Assume that a semigroup $S$ has a Brandt semigroup of left I-quotients $Q$. Put $\textbf{C}=S^*\cup E(Q)$ and  $\mathbb{G}= Q^*$. Define multiplication on $\mathbb{G}$ as  defined before Lemma~\ref{fromprimitivetog}. By Lemma~\ref{fromprimitivetog},    $\textbf{C}$ is a left order in $\mathbb{G}$ which is an inductive groupoid. We claim that $\mathbb{G}$ is connected. Suppose that $a^{-1}a$ and $b^{-1}b$ are two identity elements of $\mathbb{G}$ for some $a,b\in \textbf{C}$. Since $S$ is a left I-order in $Q$, it follows that  for $a,b\in S$ there exist $c,d \in S$ such that $ca\,\mathcal{R}^*\,d\,\lambda\,b$, by (E). As $ca\neq 0$ in $S$ we have that $ca$ is defined in $\textbf{C}$ and so $c^{-1}c=aa^{-1}$. Since $ca\,\mathcal{R}\,d$ in $Q$, by Lemma 2.4 in \cite{GG}.
 We have \[\d(ca)=\d(c)=cc^{-1}=c(c^{-1}c)c^{-1}=c(aa^{-1})c^{-1}=(ca)(ca)^{-1}=dd^{-1}=\d(d).\]
As $d\,\lambda\,b$ in $S$ so that $xd=yb\neq 0$ for some $x,y \in S$. Hence $xd$ and $yb$ are defined in $\textbf{C}$ and $xd=yb$ in $\textbf{C}$ so that $\r(b)=\r(d)$. It is clear that $\r(ca)=\r(a)$. By Lemma~\ref{connected},  $\mathbb{G}$ is connected. We have established our claim. We have

\begin{lemma}\label{fromNtoB}
  If a semigroup $S$ is a left I-order in a Brandt semigroup $Q$, then $\textbf{C}$ is a left order in $\mathbb{G}$ where $\textbf{C}=S^*\cup E(Q)$ and $\mathbb{G}=Q^*$.
 \end{lemma} 

Now, we aim to turn left orders in connected groupoids to Brandt semigroups of left I-quotients.\\
\\  
 Suppose that $\textbf{C}$ is a left order in a connected groupoid  $\mathbb{G}$. Let $S=\textbf{C}\cup \{0\}$ and $Q= \mathbb{G}\cup \{0\}$. Define multiplication on $\mathbb{G}$ as it defined before Lemma~\ref{fromgroupidtop}. By Lemma~\ref{fromgroupidtop},    $S$ is a left I-order in $Q$ where $Q$ is a primitive inverse semigroup. We claim that $S$ satisfies (E). In other words, we claim that $Q$ is a Brandt semigroup.\\
 Let $a,b\in S$ so that $a^{-1}a,b^{-1}b\in E(Q)=\mathbb{G}_0$ as $\mathbb{G}$ is connected, there are $c,d \in \textbf{\textbf{C}}$ such that $\d(c)=\d(d), \ \r(a)=\r(c)$  and $\r(b)=\r(d)$. As $\d(ca)=\d(c)=\d(d)$ in $\mathbb{G}$ we have that $ca\,\mathcal{R}\,d$ in $Q$ so that $ca\,\mathcal{R}^*\,d$ in $S$, by Lemma 2.4 in \cite{GG}. 
 Since $\textbf{C}$ is a left order in  $\mathbb{G}$ we have that $\textbf{C}$ is right reversible  and as $\r(b)=\r(d)$ we have that $xb=yd$ is defined for some $x,y \in \textbf{\textbf{C}}$ so that $xb=yd\neq 0$ in $S$ and so $d\,\lambda\,b$.  Thus we have established our claim. Hence we have
 
 \begin{lemma}\label{fromconntoBra}
 If a category $\textbf{C}$ is a left order in a connected groupoid $\mathbb{G}$, then $S=\textbf{C}\cup \{0\}$ is a left I-order in a Brandt semigroup $Q= \mathbb{G}\cup \{0\}$. \end{lemma}

\section{Inductive $\omega$-groupoids of left quotients}\label{w-ordered}

In this section we are concerned with a special class of connected groupoids of left quotients. This class is associated to a bisimple inverse $\omega$-semigroup.  We begin by describing the well known \emph{Bruck-Reilly extension}.\\
\\
Let $G$ be a group and $\theta$ be an endomorphism of $G$. The Bruck-Reilly extension $BR(G,\theta)$
of $G$ with respect to $\theta$ is the set $\mathbb{N}^0\times G\times \mathbb{N}^0$ with the binary operation:
\[(m, a, n)(p, b, q) = \big(m- n + s, (a\theta ^{s-n})(b\theta^{s-p}),q-p + s\big)\]
where $s = \max(n, p)$. The idempotents of $BR(G,\theta)$ are the elements of the form $(n, 1_G, n)$
where $n \in \mathbb{N}^0$. 

\begin{theorem}\cite{clifford} Every bisimple inverse $\omega$-semigroup is isomorphic to some
Bruck-Reilly extension of a group $G$ determined by an endomorphism of $G$.\end{theorem}

We denote $(\omega,\leq)$ the poset consisting of the natural numbers under the dual of the usual partial order.\par
From Proposition~\ref{prim2nd}, we know that any connected groupoid has the form $I\times G\times I$ where $I$ is a nonempty set and $G$ is a group. If we chose $I$ to be $\mathbb{N}^0$, then we obtain a connected groupoid $\mathcal{T}=\mathbb{N}^0 \times G \times\mathbb{N}^0$. The identities of $\mathcal{T}$ have the form $(i,1,i)$ and we have the $\omega$-ordering on the identities  as follows
\[(a,1,a)\leq (b,1,b) \Longleftrightarrow a\geq b,\] that is,
\[(0, 1, 0) \geq (1, 1, 1) > (2, 1, 2) > (3, 1, 3) > . . . .\]
It is shown in \cite{Gel}, that  $\mathcal{T}$ is inductive, and it is called an \textbf{\emph{inductive $\omega$-groupoid}}.\\
 \\
 For $(a,g,b)\in \mathcal{T}$ we have that $\d((a,g,b))=(a,1,a)$ and $\r((a,g,b))=(b,1,b)$.
 It is shown in \cite{Gel} that with a bisimple inverse $\omega$-semigroup $Q$ we can associate an inductive  $\omega$-groupoid isomomorhic to one associated to a Bruck-Reilly extension.\par
 \medskip
A natural question to ask at this point is, what are the necessary and sufficient conditions of a category to have an inductive $\omega$-groupoid of left quotients. By Lemmas~\ref{connected} and  Lemma~\ref{identities},  the first part of the following lemma is clear.
 
 \begin{lemma}\label{w-inductive}
A category $\textbf{C}$ is a left order in an inductive $\omega$-groupoid $\mathcal{T}$  if and only if  $\textbf{C}$ satisfies the following conditions:

$(i)$ $\textbf{C}$ is right reversible;

$(ii)$ $\textbf{C}$ is cancellative;

$(iii)$ $\textbf{C}$ has the connected condition;

$(iv)$ $\textbf{C}_0$ is an $\omega$-chain.
\end{lemma}
\begin{proof}
Suppose that $\textbf{C}$ satisfies the Conditions $(i)-(iv)$. By Corollary~\ref{catequ} and Lemma~\ref{connected},  we have that $\textbf{C}$ is a left order in a connected groupoid $\mathbb{G}$. By Lemma~\ref{rightreversiblecond}, $\mathbb{G}$ has the same identities as $\textbf{C}$. Hence $\mathbb{G}$ is an inductive $\omega$-groupoid, by $(iv)$.\par
Conversely, If $\textbf{C}$ is a left order in an inductive $\omega$-groupoid $\mathcal{T}$, then by Corollary~\ref{catequ},  $\textbf{C}$ is right reversible and cancellative. By Theorem~\ref{Markproof}, $\textbf{C}_0$ is an $\omega$-chain. Since $\mathcal{T}$ is connected we have that $\textbf{C}$ has the connected condition, by Lemma~\ref{connected}.
\end{proof} 
 
Following \cite{GG}, let $\mathcal{B}$ be the bicyclic monoid. Consider a semigroup $S$ together with a homomorphism $\varphi: S\longrightarrow \mathcal{B}$. We define functions $l,r:S \longrightarrow \mathbb{N}^0$ by  $a\varphi=\big(r(a),l(a)\big)$. We also put $H_{i,j}=(i,j)\varphi^{-1}$, so that $S$ is a disjoint union of subsets of the $H_{i,j}$ and \[H_{i,j}=\{a\in S: r(a)=i, \ l(a)=j\}.\]

It is well known that  $\mathcal{H}$ is a congruence on any bisimple inverse  $\omega$-semigroup  $Q$ and $Q/\mathcal{H}\cong \mathcal{B}$ where  $\mathcal{B}$ is the bicyclic semigroup. Let  $\overline{\varphi} :Q\longrightarrow \mathcal{B}$ be a surjective  homomorphism with  $Ker \overline{\varphi}=\mathcal{H}$.  As above  we will index the  $\mathcal{H}$-classes of $Q$ by putting  $H_{i,j}=(i,j)\overline{\varphi}^{-1}$.\par 
\bigskip
Let  $S$ be a left I-order in  $Q$. Let $\varphi=\overline{\varphi}|_S$ so that  $\varphi$ is a homomorphism from  $S$ to  $\mathcal{B}$. Unfortunately, this homomorphism is not surjective in general, since  $S$ need not intersect every  $\mathcal{H}$-class of  $Q$. But we can as above index the elements of  $S$.

\begin{theorem}\label{bisimplew}\cite{GG}
A semigroup  $S$  is a left I-order in a bisimple inverse $\omega$-semigroup  $Q$  if and only if  $S$  satisfies the following conditions:

$(A)$ There is a homomorphism  $\varphi:S\longrightarrow \mathcal{B}$ such that $S\varphi$ is a left I-order in $\mathcal{B}$;
  
$(B)$ For  $x   , y , a \in S$,  
\[(i) \ l(x), l(y) \geqslant r(a)\ \mbox{and} \ xa=ya \ \mbox{implies} \ x=y, \]
\[(ii) \ r(x), r(y) \geqslant l(a)\ \mbox{and} \ ax=ay \ \mbox{implies} \ x=y. \]

$(C)$ For any $b,c\in S$, there exist \ $x, y \in S$\ such that $xb = yc$ where
\[x \in H_{r(x),r(b)-l(b)+max\big(l(b), l(c)\big)},  y \in H_{r(x),r(c)-l(c)+max\big(l(b), l(c)\big)}.\]
\end{theorem}

Suppose that $S$ is a left I-order in a bisimple inverse $\omega$-semigroup  $Q$. Let $\mathcal{T}$ be an inductive $\omega$-groupoid associated to $Q$. The restricted product $\circ$ is defined on $Q$ by the rule that 
\[p\circ q = \begin{cases}
pq, & \mbox{if} \ p^{-1}p=qq^{-1}; \\
\mbox{undefined,} & \mbox{otherwise,}.  
\end{cases}\]
Then $\mathcal{T}=(Q,\circ)$ is an inductive groupoid. It is clear that $\textbf{C}=(S\cup E(Q),\circ)$ is a subcategory of $\mathcal{T}$. For if, $a\in S$ and $e\in E(Q)$ and $e\circ a (a\circ e)$ is defined in $\mathcal{T}$, then $e=aa^{-1} (a^{-1}a=e)$ so that $e\circ a=a (a\circ e=a) \in \textbf{C}$. It is clear that $\textbf{C}$  is a left order in $\mathcal{T}$. We have

\begin{lemma}
If a semigroup $S$ is a left I-order in  a bisimple inverse $\omega$-semigroup  $Q$, then $\textbf{C}=(S\cup E(Q),\circ)$ is a left order in the inductive groupoid $\mathcal{T}=(Q,\circ)$.
\end{lemma}
 
If a category $\textbf{C}$ is a left order in an inductive $\omega$-groupoid $\mathcal{T}$, then by Lemma~\ref{fromindtoinv}, $(\textbf{C},\otimes)$ is a left I-order in $(\mathcal{T},\otimes)$. Hence $(\mathcal{T},\otimes)$ is a bisimple inverse $\omega$-semigroup so that $(\textbf{C},\otimes)$ satisfies conditions (A), (B) and (C).% Moreover, $(\mathcal{T},\otimes)$ is proper **.

\begin{lemma}
If a category $\textbf{C}$ is a left order in an inductive $\omega$-groupoid $\mathcal{T}$, then the semigroup $(\textbf{C},\otimes)$ is a left I-order in $(\mathcal{T},\otimes)$.
\end{lemma}

 \section{Left q-orders in groupoids}\label{leftq-order}
 In Section~\ref{primitiveord} we showed that if a category $\textbf{C}$ is a left order in a groupoid $\mathbb{G}$, then $S=\textbf{C}^0$ is a left I-order in the primitive inverse semigroup $Q=\mathbb{G}^0$. The problem is, $S$ is a full subsemigroup of $Q$ as $\textbf{C}$ and $\mathbb{G}$ have the same set of identities. To solve this problem we need to delete the identites of $\textbf{C}$. In other words, we need to consider semigroupids.\par
 \bigskip
This section is entirely devoted to proving Theorem ~\ref{cateqqu} which gives a characterisation of semigroupoids which have a groupoid of left q-quotients.\par
 \medskip
We recall that a semigroupoid $\mathfrak{C}$ is given by
 
$(a)$ a set $\Ob (\mathfrak{C})$ of objects;

$(b)$ for each pair $(a,b)$ of objects, a set $ \mathfrak{C}(a,b)$ of homomorphisms;

$(c)$ for each triple $(a,b,c)$ of objects, a mapping from $\mathfrak{C}(a,b)\times \mathfrak{C}(b,c)$ into $\mathfrak{C}(a,c)$ which associates to each $u\in \mathfrak{C}(a,b)$ and $v\in \mathfrak{C}(b,c)$ the associative composition $uv\in \mathfrak{C}(a,c)$. We say that $u\in \mathfrak{C}(a,b)$ has \emph{domain} $\textbf{\dom}(u)=a$ and \emph{codomain} $\textbf{\cod}(u)=b$. We write $\Hom(\mathfrak{C})=\bigcup _{a,b\in {\rm\Ob}(\mathfrak{C})}\mathfrak{C}(a,b)$. \par
 \medskip
A semigroupoid is \emph{connected} if $\mathfrak{C}(u,v)\neq \emptyset$ for all $u,v \in \mathfrak{C}$. We recall that a semigroupoid is said to be \emph{right cancellative} if $\exists x\cdot a,\exists y\cdot a $ and $xa=ya$ implies $x=y$. A \emph{left cancellative semigroupoid} is defined dually. A \emph{cancellative semigroupoid} is one which is both left and right cancellative.\\
\\

It is noted in \cite{Kambites} that for a semigroupoid $\mathfrak{C}$ one can adjoint a new element $0$ not in  $\mathfrak{C}$  such that $\mathfrak{C}\cup \{0\}$ is a semigroup with multiplication 
\begin{displaymath}
ab = \left\{ \begin{array}{ll}
\mbox{the} \ \mathfrak{C}\mbox{-product} \;  ab \ \mbox{if} \; a,b \in \mathfrak{C} \;  \mbox{and} \ \cod(a)=\dom(b); \\
\quad 0, \hspace{2cm} \mbox{otherwise.}
\end{array} \right.\end{displaymath} 
Moreover, $\mathfrak{C}\cup \{0\}$ is categorical at $0$.\par
 \medskip
A semigroupoid $\mathfrak{C}$ is said to be \emph{right reversible} if for all $a,b\in \mathfrak{C}$, with $\textbf{\cod}(a)=\textbf{\cod}(b)$, there exist $p,q\in \mathfrak{C}$ such that $pa=qb$ where $\textbf{\cod}(a)$ and $\textbf{\cod}(b)$ are codomain $a$ and $b$ respectively. \par
 \medskip
From Lemmas ~\ref{rightreversiblecond} and ~\ref{equality}, we can easily deduce the following lemmas. The proofs are clear and will be omitted.
\begin{lemma}\label{mainq-order}Let  $\mathfrak{C}$ be left order in a groupoid $\mathbb{G}$. Then

$(i)$ $\mathfrak{C}$  is cancellative;

$(ii)$ $\mathfrak{C}$ is a right reversible;
\end{lemma}

\begin{lemma}\label{qequality}
Suppose that $\mathbb{G}$ is a groupoid of left quotients of $\mathfrak{C}$. Then for all $a,b,c,d\in \mathfrak{C}$ the following are equivalent.

$(i)$ \ $a^{-1}b=c^{-1}d$;

$(ii)$ \ there exist $ x,y \in \mathfrak{C}$ such that  $xa=yc$ and $xb=yd$;

$(iii)$ \  $\cod(a)=\cod(c), \cod(b)=\cod(d)$ and for all $x,y \in \mathfrak{C}$ we have  $xa=yc \Longleftrightarrow xb=yd$.
                     
\end{lemma}

We now state the main result of this section.

\begin{theorem}\label{cateqqu}
 A subsemigroupoid $\mathfrak{C}$ is  a left order in a groupoid $\mathbb{G}$ if and only if    satisfies the following conditions:
 
 $(A)$ $\mathfrak{C}$ is right reversible;
 
 $(B)$ $\mathfrak{C}$ cancellative;
 
 $(C)$ for all $a \in \mathfrak{C}$ there exists $x\in \mathfrak{C}$  such that $xa$ is defined.
\end{theorem}
\begin{proof}
If $\mathfrak{C}$ is a left order in $\mathbb{G}$, then by Lemmas~\ref{mainq-order},  $\mathfrak{C}$ is right reversible and cancellative. Hence $(A)$ and $(B)$ hold. It remains to prove $(C)$, let $a$ be any element in $\mathfrak{C}$. As  $\mathfrak{C}$ is a left q-order in $\mathbb{G}$ we have that $a=x^{-1}y$ for some $x,y \in \mathbb{G}$. Thus $xa=xx^{-1}y=y$ is defined in $\mathfrak{C}$ as required. \par
Conversely, we suppose that $\mathfrak{C}$  satisfies Conditions (A)-(C). We aim to construct a groupoid $\mathbb{G}$ in which $\mathfrak{C}$ is embedded as a left q-order in $\mathbb{G}$.  Let
\[\widetilde{\mathbb{G}}=\{(a,b)\in \mathfrak{C}\times \mathfrak{C}: \dom(a)=\dom(b)\}.\]
 Define a relation $(a,b)\sim (c,d)$ on $\widetilde{\mathbb{G}}$ by 
\[(a,b)\sim (c,d) \Longleftrightarrow \ \mbox{there exist}\ x,y \in \mathfrak{C}\ \mbox{such that} \ xa=yc\ \mbox{and}\ xb=yd.\]
Notice that if $(a,b)\sim (c,d)$, then $\cod(a)=\cod(c)$ and $\cod(b)=\cod(d)$.\par
\medskip
The following lemma is clear.
\begin{lemma}
The relation $\sim$ defined above is an equivalence relation.
\end{lemma}

Let $[a,b]$ denote the $\sim$-equivalence class of $(a,b)$. On $\mathbb{G}=\widetilde{\mathbb{G}}/$$\sim$ we define a product by

\[[a,b][c,d] = \begin{cases}
[xa,yd] & \ \mbox{if}\ \cod(b)=\cod(c) \ \mbox{and}\ xb=yc \  \mbox{for some}\ x,y \in \mathfrak{C}; \\
\mbox{undefined} & \mbox{otherwise},  
\end{cases}\] 
and so we have the following diagram
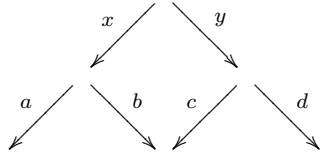
\[\begin{xy}
\xymatrix{
&&  \ar[dr]^y \ar[dl]_x & \\
& \ar[dl]_a\ar[dr]^b & & \ar[dl]_c\ar[dr]^d&&
 \\  & &  & & }
\end{xy}
\]

The proofs of the following Lemmas are similar to the proofs of Lemmas~\ref{welldefinedgroupoid} and~\ref{associative}.
\begin{lemma} 
The multiplication is well-defined.
\end{lemma}

\begin{lemma}\label{qassociative}
The multiplication is associative.
\end{lemma}

Before proceeding with the proof we insert a lemma which we shall use often.

\begin{lemma}\label{supportive} If $[a,b] \in \mathbb{G}$, then $[xa,xb]=[a,b]$ for some $x\in \mathfrak{C}$ such that  $xa$ is defined in $\mathfrak{C}$.\end{lemma}
\begin{proof}
Let $[a,b] \in \mathbb{G}$ where $xa$ is defined in $\mathfrak{C}$ for some $x\in \mathfrak{C}$. It is clear that $[xa,xb]\in \mathbb{G}$. Since $\mathfrak{C}$ is right reversible and $\cod(xa)=\cod(a)$ we have that $txa=ra$ for some $t,r \in \mathfrak{C}$, by the right reversiblility of $\mathfrak{C}$. By cancellativity in $\mathfrak{C}$ we get $tx=r$. Since $\dom(a)=\dom(b)$, it follows that $rb$ is defined. Hence   $txb=rb$. Thus $[a,b]=[xa,xb]$. \end{proof}

From Lemma~\ref{identities}, we can easily deduce the following lemma.

\begin{lemma}\label{identitiessemi}
 The identities of $\mathbb{G}$ have the form $[a,a]$ where $a\in \mathfrak{C}$.\end{lemma}

 Suppose that $[a,b]\in \mathbb{G}$. By Condition (C), there exists $x\in \mathfrak{C}$ such that $xa$ is defined. Using Lemma~\ref{supportive} and the definition of multiplication, we get 
 \[[a,a][a,b]=[xa,xb]=[a,b].\]
Similarly, 
\[[a,b][b,b]=[ya,yb]=[a,b],\]
for some $y \in \mathfrak{C}$ such that $yb$ is defined. Hence $\d([a,b])=[a,a]$ and $\r([a,b])=[b,b]$.\par
\medskip
 By the above argument and Lemma~\ref{qassociative}, the following lemma is clear.
 \begin{lemma} $\mathbb{G}$ is a category.\end{lemma}

If $[a,b]\in \mathbb{G}$, then it is clear that $[b,a] \in \mathbb{G}$. By Condition (C) and Lemma~\ref{supportive},  
 \[[a,b][b,a]=[ya,ya]=[a,a]=\d([a,b])\]
 for some $y \in \mathfrak{C}$ such that $yb$ is defined. Similarly, $[b,a][a,b]=[b,b]=\r([a,b])$. That is, $[b,a]$ is the inverse of $[a,b]$ in $\mathbb{G}$. Thus we have

  \begin{lemma} $\mathbb{G}$ is a groupoid.\end{lemma}

If  $a \in \mathfrak{C}$, by (C) there exists $x\in \mathfrak{C}$ such that  $xa$ is defined. It is clear that   $[x,xa] \in \mathbb{G}$. If  $(y,ya) \in \widetilde{\mathbb{G}}$, then as $\cod(xa)=\cod(ya)$   there exist  $u,v \in \mathfrak{C}$ such that   $uxa=vya$. Since $\mathfrak{C}$ is cancellative  we have  that $ux=vy$, that is,  $[x,xa]=[y,ya]$. Thus we have a well defined map $\theta: \mathfrak{C} \longrightarrow \mathbb{G}$ given by $a\theta=[x,xa]$ where $x \in \mathfrak{C}$ such that $xa$ is defined.

 \begin{lemma}\label{qembedding} The mapping $\theta$ is an embedding of $\mathfrak{C}$ in $\mathbb{G}$.
 \end{lemma} 
\begin{proof}
 To show that $\theta$ is one-to-one, let $[x,xa]=[y,yb]$ for some $x,y \in \mathfrak{C}$ so that $ux=vy$ and $uxa=vyb$ for some $u,v\in \mathfrak{C}$. Hence $a=b$.\\
\\
For $a,b\in \mathfrak{C}$ such that $ab$ is defined. We have 
\[\begin{array}{rcl}a\theta b\theta &=& [x,xa][y,yb]\\ &=&[ux,vyb] \\ &=&[ux,uxab]\\ &=&(ab)\theta,\end{array}\]
where $uxa=vy$ for some $u,v \in \mathfrak{C}$. Thus $\theta$ is a homomorphism.
\end{proof}

By Lemma~\ref{qembedding}, we can regard $\mathfrak{C}$ as a subsemigroupoid  of $\mathbb{G}$. Let $[a,b]\in \mathbb{G}$ and  $a\theta=[x,xa]$, $b\theta=[y,yb]$ where $xa$ and $yb$ are defined for some $x,y\in \mathfrak{C}$. As $\cod(x)=\dom(a)=\dom(b)=\cod(y)$, by the right reversibility of $\mathfrak{C}$ there are elements   $u,v\in \mathfrak{C}$ with $ux=vy$. Hence
 \[\begin{array}{rcl}(a\theta)^{-1} (b\theta) &=& [x,xa]^{-1}[y,yb]\\ &=& [xa,x][y,yb]\\ &=& [uxa,vyb]\\ &=&[uxa,uxb]\hspace{1.3cm} \mbox{by Lemma~\ref{supportive}},\\ &=&[a,b] .\end{array}\]
Hence $\mathfrak{C}$ is a left order in  $\mathbb{G}$. This completes the proof of Theorem~\ref{cateqqu}.\end{proof}

The following lemma can be deduced from Corollary~\ref{coruniqu}.
 \begin{lemma}
If a semigroupoid $\mathfrak{C}$ is a left order in groupoids $\mathbb{G}$ and  $\mathbb{P}$,
then $\mathbb{G}$ and  $\mathbb{P}$ are isomorphic by an isomorphism which restricts to the identity map on $\textbf{C}$.
 \end{lemma}

Let $\mathfrak{C}$ be a semigroupoid. We say that it satisfies the \textit{connected condition} if,  for any $a,b\in \mathfrak{C}$ there exist $c,d\in \mathfrak{C}$ with $\dom(c)=\dom(d)$ such that  $\cod(c)=\cod(a)$  and  $\cod(d)=\cod(b)$.\\
\\

We conclude this section with the following results. The proofs can be deduced easily from those for categories in Sections ~\ref{primitiveord} and   ~\ref{connectedored}.
\begin{lemma}\label{fromprimitivetogsemi}
Let $S$ be a semigroup. If $S$ is a left I-order in a  primitive inverse (Brandt) semigroup $Q$, then $\mathfrak{C}=S^*$ is a left q-order in the (connected) groupoid $Q^*$.
 \end{lemma}

\begin{lemma}\label{fromgroupidtopsemi}
Let $\mathfrak{C}$ be a semigroupoid. If $\mathfrak{C}$ is a left q-order in a (connected) groupoid $\mathbb{G}$, then $S=\mathfrak{C}\cup \{0\}$ is a left I-order in (Brandt) primitive inverse semigroup $Q=\mathbb{G} \cup\{0\}$.
 \end{lemma}

\end{document}